\pdfoutput=1
\RequirePackage{ifpdf}
\ifpdf 
\documentclass[pdftex]{sigma}
\else
\documentclass{sigma}
\fi

\usepackage{enumerate}

\numberwithin{equation}{section}

\newtheorem{thm}{Theorem}[section]

\newtheorem{lem}[thm]{Lemma}
\newtheorem{prop}[thm]{Proposition}
\theoremstyle{definition}

\newtheorem{Remark}[thm]{Remark}

\begin{document}

\newcommand{\arXivNumber}{1502.05252}

\allowdisplaybreaks

\renewcommand{\PaperNumber}{054}

\FirstPageHeading

\ShortArticleName{Eigenvalue Estimates of the ${\mathop{\rm spin}^c}$ Dirac Operator and Harmonic Forms}

\ArticleName{Eigenvalue Estimates of the $\boldsymbol{{\mathop{\rm spin}^c}}$ Dirac Operator\\ and Harmonic Forms on K\"ahler--Einstein Manifolds}

\Author{Roger NAKAD~$^\dag$ and Mihaela PILCA~$^{\ddag\S}$}

\AuthorNameForHeading{R.~Nakad and M.~Pilca}

\Address{$^\dag$~Notre Dame University-Louaiz\'e, Faculty of Natural and Applied Sciences,\\
 \hphantom{$^\dag$}~Department of Mathematics and Statistics,  P.O. Box 72, Zouk Mikael, Lebanon}
\EmailD{\href{mailto:rnakad@ndu.edu.lb}{rnakad@ndu.edu.lb}}
\URLaddressD{\url{http://www.iecn.u-nancy.fr/~nakad/}}

\Address{$^\ddag$~Fakult\"at f\"ur Mathematik,  Universit\"at Regensburg,\\
\hphantom{$^\ddag$}~Universit\"atsstra{\ss}e~31, 93040 Regensburg, Germany}
\EmailD{\href{mailto:Mihaela.Pilca@mathematik.uni-regensburg.de}{Mihaela.Pilca@mathematik.uni-regensburg.de}}
\URLaddressD{\url{http://www.mathematik.uni-regensburg.de/pilca/}}

\Address{$^\S$~Institute of Mathematics ``Simion Stoilow'' of the Romanian Academy,\\
\hphantom{$^\S$}~21, Calea Grivitei Str, 010702-Bucharest, Romania}

\ArticleDates{Received March 03, 2015, in f\/inal form July 02, 2015; Published online July 14, 2015}

\Abstract{We establish a lower bound for the eigenvalues of the Dirac operator def\/ined on a~compact K\"ahler--Einstein manifold of positive scalar curvature and endowed with particular ${\mathop{\rm spin}^c}$ structures. The limiting case is characteri\-zed  by the existence of K\"ahlerian Killing ${\mathop{\rm spin}^c}$ spinors in a certain subbundle of the spinor bundle. Moreover, we show that the Clif\/ford multiplication between an ef\/fective harmonic form and a K\"ahlerian Killing ${\mathop{\rm spin}^c}$ spinor f\/ield vanishes.
This extends to the ${\mathop{\rm spin}^c}$ case the result of  A.~Moroianu stating that, on a~compact K\"ahler--Einstein manifold of complex dimension $4\ell+3$ carrying a complex contact structure, the Clif\/ford multiplication between an ef\/fective harmonic form and a~K\"ahlerian Killing spinor is zero.}


\Classification{53C27;  53C25; 53C55; 58J50; 83C60}

\section{Introduction}

The geometry and topology of a compact Riemannian spin manifold $(M^n, g)$ are strongly related to the existence of special spinor f\/ields and thus, to the spectral properties of a fundamental operator called the Dirac operator $D$ \cite{AS, lich0}.
A.~Lichnerowicz \cite{lich0} proved, under the weak condition of the positivity of the scalar curvature, that the kernel of the Dirac operator is trivial. Th.~Friedrich \cite{fr1} gave the following lower bound for the f\/irst eigenvalue $\lambda$ of $D$ on a~compact Riemannian spin manifold $(M^n, g)$:
\begin{gather}\label{fried}
\lambda^2\geq\frac{n}{4(n-1)}\underset{M}{\inf}\, S,
\end{gather}
where $S$ denotes the  scalar curvature, assumed to be nonnegative. Equality holds if and only if the corresponding eigenspinor $\varphi$ is parallel (if $\lambda=0$) or a Killing spinor of Killing constant $-\frac{\lambda}{n}$ (if $\lambda\neq0$), i.e., if $\nabla_X\varphi = -\frac{\lambda}{n} X\cdot\varphi$, for all vector f\/ields $X$, where ``$\cdot$'' denotes the Clif\/ford multiplication and $\nabla$  is the spinorial Levi-Civita connection on the spinor bundle $\Sigma M$ (see also~\cite{hijconf}). Killing (resp.\ parallel) spinors force the underlying metric to be Einstein  (resp.\ Ricci f\/lat).  The classif\/ication of complete simply-connected Riemannian spin manifolds with real Killing (resp.\ parallel) spinors was done by C.~B\"ar~\cite{baer} (resp.\ M.Y.~Wang~\cite{wang}). Useful geometric information has been also obtained by restricting parallel and Killing spinors to hypersurfaces \cite{Ba, HMZ1, HMZ2, HMZ02, HZ1, HZ2}. O.~Hijazi  proved  that the Clif\/ford multiplication between  a harmonic  $k$-form $\beta$  ($k \neq 0, n$) and a Killing spinor vanishes. In particular, the equality case  in~\eqref{fried} cannot be attained on a  K\"ahler spin manifold, since the Clif\/ford multiplication between the K\"ahler form and a Killing spinor is never zero.
Indeed, on a K\"ahler compact manifold $(M^{2m},g,J)$ of complex dimension $m$ and  complex structure $J$, K.-D.~Kirchberg~\cite{kirch86}  showed that the f\/irst  eigenvalue $\lambda$ of the Dirac operator satisf\/ies
\begin{gather}\label{kirchoddeven}
\lambda^2\geq\begin{cases}
\dfrac{m+1}{4m}\underset{M}{\inf}\,  S,&  \text{if $m$ is odd,} \vspace{1mm}\\
\dfrac{m}{4(m-1)}\underset{M}{\inf} \, S,&   \text{if $m$ is even.}
\end{cases}
\end{gather}
Kirchberg's estimates rely essentially on the decomposition of $\Sigma M$ under the action of the K\"ahler form $\Omega$. In fact, we have $\Sigma M = \oplus_{r=0}^m\Sigma_r M$, where~$\Sigma_r M$ is the eigenbundle corresponding to the eigenvalue $i(2r-m)$ of $\Omega$. The limiting manifolds of (\ref{kirchoddeven}) are also characterized by the existence of
spinors satisfying a certain dif\/ferential equation similar to the one fulf\/illed by Killing spinors.
More precisely, in odd complex dimension $m=2\ell+1$, it is proved in \cite{hij, kirch2,kirch} that the metric is Einstein and the corresponding eigenspinor $\varphi$ of $\lambda$ is a K\"ahlerian Killing spinor,  i.e., $\varphi=\varphi_{\ell}+\varphi_{\ell+1}\in\Gamma(\Sigma_{\ell}  M\oplus \Sigma_{\ell+1}  M)$ and
it satisf\/ies
\begin{gather}
\begin{split}
&\nabla_X \varphi_\ell = -\frac{\lambda}{2(m+1)} (X+iJX)\cdot \varphi_{\ell+1}, \\
&\nabla_X \varphi_{\ell+1} = -\frac{\lambda}{2(m+1)} (X-iJX)\cdot \varphi_{\ell},
\end{split}\label{ecodd}
\end{gather}
for any vector f\/ield $X$.  We point out that the existence of spinors of the form $\varphi=\varphi_{\ell'}+\varphi_{\ell'+1}\in\Gamma(\Sigma_{\ell'}M\oplus\Sigma_{\ell'+1}M)$
satisfying \eqref{ecodd}, implies that $m$ is odd and they lie in the middle,  i.e., $l' =\frac{m-1}{2}$. If the complex dimension is even, $m=2\ell$, the limiting manifolds are characterized by constant scalar curvature and the existence of so-called anti-holomorphic K\"ahlerian twistor spinors $\varphi_{\ell-1}\in\Gamma(\Sigma_{\ell-1}M)$, i.e., satisfying for any vector f\/ield $X$: $\nabla_X \varphi_{\ell-1}= -\frac{1}{2m}(X+iJX)\cdot D\varphi_{\ell-1}$.
The limiting manifolds for Kirchberg's inequalities \eqref{kirchoddeven} have been geometrically described by
A.~Moroianu in~\cite{am_odd} for~$m$ odd and in~\cite{am_even} for~$m$ even. In \cite{pilcapaper}, this result is extended to limiting manifolds of the so-called ref\/ined Kirchberg inequalities, obtained by restricting the square of the Dirac operator to the eigenbundles~$\Sigma_r M$. When~$m$ is even, the limiting manifold cannot be Einstein. Thus,  on compact K\"ahler--Einstein manifolds of even complex dimension, K.-D.~Kirchberg~\cite{kircheven} impro\-ved~\eqref{kirchoddeven}  to the following lower bound
\begin{gather}\label{kirchke}
\lambda^2 \geq \frac{m+2}{4m} S.
\end{gather}
Equality is characterized by the existence of holomorphic or anti-holomorphic spinors.
When $m$ is odd, A.~Moroianu extended the above mentioned result of O.~Hijazi to K\"ahler manifolds, by showing that
the Clif\/ford multiplication between a harmonic  ef\/fective form of nonzero degree and a K\"ahlerian Killing spinor vanishes.
We recall that the manifolds  of complex dimension $m =4\ell+3$ admitting K\"ahlerian Killing spinors are exactly
the K\"ahler--Einstein manifolds carrying a complex contact structure (cf.~\cite{ks, am_odd, ms}).

In the present paper, we extend this result of A.~Moroianu to K\"ahlerian Killing ${\mathop{\rm spin}^c}$ spinors (see Theorem~\ref{eff}). In this more general setting dif\/f\/iculties occur due to the fact that the connection on the ${\mathop{\rm spin}^c}$ bundle, hence its curvature, the Dirac operator and its spectrum,  do not only depend on the geometry of the manifold, but also on the connection
of the auxiliary line bundle associated with the ${\mathop{\rm spin}^c}$ structure.

\looseness=-1
$ \mathrm{Spin}^c$ geometry became an active f\/ield of research with the advent
of Seiberg--Witten theory, which has many applications to $4$-dimensional geometry and topology \cite{don,gursky,lebrun1,lebrun2,SW3, SW2}.
From an intrinsic point of view,
almost complex, Sasaki and some classes of CR manifolds carry a~canonical ${\mathop{\rm spin}^c}$ structure.
In particular, every K\"ahler manifold is ${\mathop{\rm spin}^c}$ but not necessarily spin.
For example, the complex projective space $\mathbb C P^m$ is spin if and only if $m$ is odd. Moreover, from the extrinsic point of view, it seems that it is more natural to work with ${\mathop{\rm spin}^c}$ structures rather than spin structu\-res~\cite{hmu, nakadthesis,JRRN}. For instance, on K\"ahler--Einstein manifolds of positive scalar curvature, O.~Hijazi, S.~Montiel and F.~Urbano \cite{hmu} constructed
${\mathop{\rm spin}^c}$ structures carrying  K\"ahlerian Killing  ${\mathop{\rm spin}^c}$ spinors, i.e., spinors satisfying \eqref{ecodd}, where the covariant derivative is the ${\mathop{\rm spin}^c}$ one. In \cite{her}, M.~Herzlich and A.~Moroianu  extended Friedrich's estimate \eqref{fried} to compact Riemannian ${\mathop{\rm spin}^c}$ manifolds. This new lower bound involves only the conformal geometry of the manifold and the curvature of the auxiliary line bundle associated with the ${\mathop{\rm spin}^c}$ structure. The limiting case is characterized by the existence of a ${\mathop{\rm spin}^c}$ Killing  or parallel spinor, such that the Clif\/ford multiplication of the curvature form of the auxiliary line bundle with this spinor is proportional to it.

In this paper, we give an  estimate for the eigenvalues of the ${\mathop{\rm spin}^c}$ Dirac operator, by restric\-ting ourselves to
compact K\"ahler--Einstein manifolds endowed with particular ${\mathop{\rm spin}^c}$ structures.
More precisely, we consider $(M^{2m},g,J)$  a compact K\"ahler--Einstein manifold of positive scalar curvature $S$ and of index $p\in\mathbb{N}^*$.  We endow $M$ with the ${\mathop{\rm spin}^c}$ structure whose auxiliary line bundle is a tensorial power $\mathcal{L}^q$ of the $p$-$th$ root $\mathcal L$ of the canonical bundle $K_M$ of $M$, where $q \in \mathbb Z$, $p+q\in2\mathbb{Z}$ and $|q|\leq p$. Our main result is the following:

\begin{thm}\label{globestimke}
 Let $(M^{2m}, g)$ be a compact K\"ahler--Einstein manifold of index $p$ and positive scalar curvature $S$, carrying the ${\mathop{\rm spin}^c}$ structure given by $\mathcal{L}^q$ with $q+p\in 2\mathbb{Z}$, where $\mathcal{L}^p=K_M$. We assume that $p \geq |q|$ and the metric is normalized such that its scalar curvature equals $4m(m+1)$.
 Then, any eigenvalue $\lambda$ of $D^2$ is bounded from below as follows
 \begin{gather}\label{global}
\lambda\geq \left(1-\frac{q^2}{p^2}\right) (m+1)^2.
 \end{gather}
Equality is attained if and only if $b:=\frac{q}{p}\cdot\frac{m+1}{2}+\frac{m-1}{2}\in\mathbb{N}$ and there exists a  K\"ahlerian Killing ${\mathop{\rm spin}^c}$ spinor in $\Gamma(\Sigma_{b}M\oplus \Sigma_{b+1} M)$.
\end{thm}

Indeed, this is a consequence of more ref\/ined estimates for the eigenvalues of the square of the ${\mathop{\rm spin}^c}$ Dirac operator restricted to
the eigenbundles~$\Sigma_r M$ of the spinor bundle (see Theorem~\ref{estimke}).
The proof of this result is based on a ref\/ined Schr\"odinger--Lichnerowicz~${\mathop{\rm spin}^c}$ formula (see Lemma~\ref{refinedd}) written
on each such eigenbundle~$\Sigma_rM$, which uses  the decomposition of the covariant derivative acting on spinors into its holomorphic  and antiholomorphic part. This formula has already been used in literature, for instance by K.-D.~Kirchberg~\cite{kircheven}. The limiting manifolds of~\eqref{global} are characterized by the existence of  K\"ahlerian Killing~${\mathop{\rm spin}^c}$ spinors in a certain subbundle~$\Sigma_r M$. In particular, this gives a positive answer to the conjectured relationship between ${\mathop{\rm spin}^c}$ K\"ahlerian Killing spinors and a lower bound
for the eigenvalues of the~${\mathop{\rm spin}^c}$ Dirac operator, as stated in \cite[Remark~16]{hmu}.

Let us mention here that the Einstein condition in Theorem~\ref{globestimke} is important in order to establish the estimate~\eqref{global}, since otherwise there is no control over the estimate of the term given by the Clif\/ford action of the curvature form of the auxiliary line bundle of the~${\mathop{\rm spin}^c}$ structure (see~\eqref{fa-ke}).

\section{Preliminaries and notation}
In this section, we set the notation and brief\/ly review  some basic facts about ${\mathop{\rm spin}^c}$ and K\"ahler geometries.
For more details we refer to the books \cite{bookspin, fr_book,spin,am_lectures}.

Let $(M^{n}, g)$ be an $n$-dimensional closed Riemannian ${\mathop{\rm spin}^c}$  manifold and denote by $\Sigma M$ its complex spinor bundle,
which has complex rank equal to $2^{[\frac{n}{2}]}$. The bundle $\Sigma M$ is endowed with a Clif\/ford multiplication denoted by ``$\cdot$'' and a scalar product denoted by $\langle \cdot, \cdot\rangle$.
Given a~${\mathop{\rm spin}^c}$ structure on~$(M^{n}, g)$, one can check that the
determinant line bundle $\mathrm{det}(\Sigma M)$ has a root~$L$ of index~$2^{[\frac
{n}{2}]-1}$. This line bundle~$L$ over~$M$ is called the auxiliary line bundle associated with the
${\mathop{\rm spin}^c}$ structure.
The connection~$\nabla^A$ on~$\Sigma M$ is the twisted connection of the one on the spinor bundle (induced
by the Levi-Civita connection) and a f\/ixed connection~$A$ on~$L$.
The ${\mathop{\rm spin}^c}$ Dirac operator $D^A$ acting on the space of sections of $\Sigma M$
is def\/ined by the composition of the connection $\nabla^A$ with the Clif\/ford
multiplication. For simplicity, we will denote $\nabla^A$ by $\nabla$ and~$D^A$ by~$D$.
In local coordinates:
\begin{gather*}
D =\sum_{j=1}^{n} e_j \cdot \nabla_{e_j},
\end{gather*}
where $\{e_j\}_{j=1,\dots, n}$ is a local orthonormal basis of~$TM$.
$D$ is a f\/irst-order elliptic operator and is formally self-adjoint with respect to the $L^2$-scalar product.
A useful tool when examining the ${\mathop{\rm spin}^c}$ Dirac operator is the
Schr\"{o}dinger--Lichnerowicz formula
\begin{gather}
D^2 = \nabla^*\nabla + \frac 14 S +
\frac{1}{2}F_A \cdot,
\label{sl}
\end{gather}
where $\nabla^*$ is the adjoint of $\nabla$ with respect to the $L^2$-scalar
product  and $F_A$ is the curvature (imaginary-valued) $2$-form on $M$ associated to the connection $A$ def\/ined on the auxiliary line bundle $L$, which acts on spinors by the extension of the Clif\/ford multiplication to dif\/ferential forms.

We recall that the complex volume element $\omega_{\mathbb{C}}=i^{[\frac{n+1}{2}]} e_1\wedge \cdots \wedge e_n$ acts as the identity on the spinor bundle if~$n$ is odd. If~$n$ is even,  $\omega_{\mathbb C}^2=1$. Thus, under the action of the complex volume element, the spinor bundle  decomposes into the  eigenspaces $\Sigma^{\pm} M$ corresponding to the~$\pm 1$ eigenspaces, the {\it positive} (resp. {\it negative}) spinors.

Every spin manifold has a trivial ${\mathop{\rm spin}^c}$ structure, by
choosing the trivial line bundle with the trivial connection whose curvature~$F_A$
vanishes.  Every K\"ahler manifold $(M^{2m},g,J)$ has a canonical
${\mathop{\rm spin}^c}$ structure induced by the complex structure~$J$.
The complexif\/ied tangent bundle decomposes into
$T^{\mathbb{C}} M = T_{1,0} M\oplus T_{0,1} M,$
 the $i$-eigenbundle (resp.~$(-i)$-eigenbundle) of the complex linear extension of~$J$.
For any vector f\/ield $X$, we denote by $X^{\pm}:=\frac{1}{2}(X\mp iJX)$ its component in $T_{1,0} M$, resp.~$T_{0,1} M$.
The spinor bundle of the canonical ${\mathop{\rm spin}^c}$ structure is def\/ined~by
\begin{gather*}
\Sigma M = \Lambda^{0,*} M =\overset{m}{\underset{r=0}{\oplus}} \Lambda^r (T_{0,1}^* M),
\end{gather*}
and its auxiliary line bundle is  $L = (K_M)^{-1}= \Lambda^m (T_{0,1}^* M)$, where $K_M=\Lambda^{m,0}M$ is the canonical bundle of $M$.
The line bundle $L$ has a canonical holomorphic connection, whose curvature form is given by $- i\rho$,
where $\rho$ is the Ricci form def\/ined, for all vector f\/ields~$X$ and~$Y$, by $\rho(X, Y) = \mathrm{Ric}(JX, Y)$ and $\mathrm{Ric}$ denotes the Ricci tensor. Let us mention here the sign convention we use to def\/ine the Riemann curvature tensor, respectively the Ricci tensor: $R_{X,Y}:=\nabla_{X}\nabla_{Y}
-\nabla_{Y}\nabla_{X}-\nabla_{[X,Y]}$ and ${\mathop{\rm Ric}}(X,Y):=\sum\limits_{j=1}^{2m}R(e_j,X,Y,e_j)$, for all vector f\/ields $X$, $Y$ on~$M$, where $\{e_j\}_{j=1,\dots, 2m}$ is a local orthonormal basis of the tangent bundle.
Similarly, one def\/ines  the so called anti-canonical ${\mathop{\rm spin}^c}$ structure, whose spinor bundle is given by
$\Lambda^{*, 0} M =\oplus_{r=0}^m \Lambda^r (T_{1, 0}^* M)$ and the auxiliary line bundle by~$K_M$. The spinor bundle of any other ${\mathop{\rm spin}^c}$ structure on $M$ can be written as
\begin{gather*}
\Sigma M =  \Lambda^{0, *} M \otimes \mathbb L,
\end{gather*}
where $\mathbb L^2 = K_M\otimes L$ and  $L$ is the auxiliary line bundle associated with this ${\mathop{\rm spin}^c}$ structure.
The K\"ahler form $\Omega$, def\/ined as $\Omega(X,Y)=g(JX,Y)$, acts on $\Sigma M$ via Clif\/ford multiplication
and this action is locally  given by
\begin{gather}\label{defomega}
\Omega\cdot \psi = \frac{1}{2} \sum_{j=1}^{2m} e_j\cdot Je_j\cdot\psi,
\end{gather}
for all $\psi\in\Gamma(\Sigma M)$, where $\{e_1, \dots, e_{2m}\}$ is a local orthonormal basis of $\mathrm{TM}$. Under this action, the spinor bundle decomposes as follows
\begin{gather}\label{decomp}
\Sigma M =\overset{m}{\underset{r=0}{\oplus}} \Sigma_r M,
\end{gather}
where $\Sigma_r M$ denotes the eigenbundle to the eigenvalue $i(2r-m)$ of $\Omega$, of complex rank $\binom{m}{k}$.
It is easy to see that $\Sigma_r M \subset \Sigma^+ M$ (resp.\ $\Sigma_r M \subset \Sigma^-M$) if and only if $r$ is even (resp.~$r$ is odd).
Moreover, for any $X \in \Gamma(TM)$ and $\varphi \in \Gamma(\Sigma_r M)$, we have $X^+ \cdot\varphi \in \Gamma(\Sigma_{r+1}M)$
and $X^-\cdot \varphi \in \Gamma(\Sigma_{r-1} M)$, with the convention $\Sigma_{-1}M=\Sigma_{m+1}M=M\times\{0\}$.
Thus, for any ${\mathop{\rm spin}^c}$ structure, we have
$\Sigma_r M = \Lambda^{0, r} M\otimes \Sigma_0 M$.
Hence, $(\Sigma_0M)^2 = K_M \otimes L,$ where $L$ is the auxiliary line bundle associated with the ${\mathop{\rm spin}^c}$ structure. For example, when the manifold is spin, we have $(\Sigma_0M)^2 = K_M$ \cite{hit, kirch86}.  For the canonical ${\mathop{\rm spin}^c}$ structure, since $L = (K_M)^{-1}$, it follows that $\Sigma_0 M$ is trivial. This yields the existence of parallel spinors (the constant functions) lying in $\Sigma_0M$, cf.~\cite{Moro1}.

Associated to the complex structure $J$, one def\/ines the following operators
\begin{gather*}
D^+ =\sum_{j=1}^{2m}e_j^+\cdot \nabla_{e_j^-},\qquad D^- =\sum_{j=1}^{2m}e_j^-\cdot \nabla_{e_j^+},
\end{gather*}
which satisfy the relations
\begin{gather*} 
D=D^+ +D^-, \qquad (D^+)^2=0, \qquad (D^-)^2=0, \qquad D^+D^- + D^-D^+ =D^2.
\end{gather*}

When restricting the Dirac operator to $\Sigma_{r}M$, it acts as
\begin{gather*}
D=D^+ +D^-\colon \ \Gamma(\Sigma_{r}M) \to \Gamma(\Sigma_{r-1}M\oplus\Sigma_{r+1}M).
\end{gather*}

Corresponding to the decomposition $TM\otimes\Sigma_{r}M\cong \Sigma_{r-1}M\oplus \Sigma_{r+1}M\oplus \mathrm{Ker}_r$,
where $\mathrm{Ker}_r$ denotes the kernel of the Clif\/ford multiplication by tangent vectors restricted to $\Sigma_{r}M$, we have, as in the spin case (for details see, e.g., \cite[equation~(2.7)]{pilcapaper}), the following Weitzenb\"ock formula relating the dif\/ferential operators acting on sections of $\Sigma_{r}M$:
\begin{gather*}
\nabla^{*}\nabla=\frac{1}{2(r+1)}D^-D^+ +\frac{1}{2(m-r+1)}D^+D^- + T_r^*T_r,
\end{gather*}
where $T_r$ is the so-called K\"ahlerian twistor operator and is def\/ined by
\begin{gather*}
T_r\varphi:= \nabla \varphi +\frac{1}{2(m-r+1)}e_j\otimes e_j^+\cdot D^-\varphi+\frac{1}{2(r+1)}e_j\otimes e_j^-\cdot D^+\varphi.
\end{gather*}
This decomposition further implies the following identity for $\varphi\in\Gamma(\Sigma_{r}M)$, by the same argument as in \cite[Lemma 2.5]{pilcapaper},
\begin{gather}\label{ident}
|\nabla\varphi|^2= \frac{1}{2(r+1)}|D^+\varphi|^2 +\frac{1}{2(m-r+1)}|D^-\varphi|^2+|T_r\varphi|^2.
\end{gather}
Hence, we have the inequality
\begin{gather}\label{ineg}
|\nabla\varphi|^2\geq \frac{1}{2(r+1)}|D^+\varphi|^2 +\frac{1}{2(m-r+1)}|D^-\varphi|^2.
\end{gather}
Equality in \eqref{ineg} is attained if and only if $T_r\varphi=0$, in which case $\varphi$ is called a K\"ahlerian twistor spinor. The Lichnerowicz--Schr\"odinger formula~\eqref{sl} yields the following:
\begin{lem}\label{alg1}
Let $(M^{2m},g,J)$ be a compact K\"ahler manifold endowed with any ${\mathop{\rm spin}^c}$ structure.  If $\varphi$ is an eigenspinor of $D^2$ with eigenvalue $\lambda$, $D^2\varphi=\lambda\varphi$, and satisfies
\begin{gather}\label{inegalg1}
|\nabla\varphi|^2\geq\frac{1}{j}|D\varphi|^2,
\end{gather}
for some real number $j>1$, and  $(S +2 F_A) \cdot \varphi = c \varphi$,
where $c$ is a positive function, then
\begin{gather}\label{firstineq1}
\lambda\geq\frac{j}{4(j-1)}\underset{M}{\inf} \, c.
\end{gather}
Moreover, equality in \eqref{firstineq1}  holds if and only if the function $c$ is constant and equality in \eqref{inegalg1} holds at all points of the manifold.
\end{lem}
Let $\{e_1,\dots, e_{2m}\}$ be  a local orthonormal basis of $M^{2m}$. We implicitly use the Einstein summation convention over repeated indices. We have the  following formulas for contractions that
hold as endomorphisms of $\Sigma_r M$:
\begin{gather}\label{kecontr1}
e_j^+\cdot e_j^-=-2r, \qquad e_j^-\cdot e_j^+=-2(m-r),
\\
\label{kecontr3}
e_j\cdot \mathrm{Ric}(e_j)=-S, \qquad e_j^-\cdot \mathrm{Ric}(e_j^+)=-\frac{S}{2}-i\rho, \qquad e_j^+\cdot \mathrm{Ric}(e_j^-)=-\frac{S}{2}+i\rho.
\end{gather}
The \looseness=-1 identities \eqref{kecontr1} follow directly from~\eqref{defomega}, which gives the action of the K\"ahler form and has~$\Sigma_r M $ as eigenspace to the eigenvalue $i(2r-m)$, implying that $ie_j\cdot Je_j=2i\Omega=-2(2r-m)$, and from the fact that $e_j\cdot e_j=-2m$. The identities~\eqref{kecontr3} are obtained from the following identities
\begin{gather*}
e_j\cdot {\mathop{\rm Ric}}(e_j)=e_j\wedge {\mathop{\rm Ric}}(e_j)-g({\mathop{\rm Ric}}(e_j),e_j)=-S,\\
 ie_j\cdot {\mathop{\rm Ric}}(Je_j)=ie_j\wedge {\mathop{\rm Ric}}(Je_j)-ig({\mathop{\rm Ric}}(Je_j),e_j)=2i\rho.
\end{gather*}
The ${\mathop{\rm spin}^c}$ Ricci identity, for any spinor $\varphi$ and any vector f\/ield $X$, is given by
\begin{gather}\label{ricident}
e_i \cdot\mathcal{R}^A_{e_i,X} \varphi= \frac{1}{2}\mathrm{Ric}(X) \cdot\varphi-\frac{1}{2} (X\lrcorner F_A)\cdot\varphi,
\end{gather}
where $\mathcal{R}^A$ denotes the ${\mathop{\rm spin}^c}$ spinorial curvature, def\/ined with the same sign convention as above, namely $\mathcal{R}^A_{X,Y}:=\nabla^A_{X}\nabla^A_{Y}
-\nabla^A_{Y}\nabla^A_{X}-\nabla^A_{[X,Y]}$. For a~proof of the ${\mathop{\rm spin}^c}$ Ricci identity we refer to \cite[Section~3.1]{fr_book}.
For any vector f\/ield $X$ parallel at the point where the computation is done, the following commutator rules hold
 \begin{gather}\label{nablax}
[\nabla_X,D]=-\frac{1}{2}\mathrm{Ric}(X)\cdot+\frac{1}{2}(X\lrcorner F_A)\cdot,
\\\label{nablax+}
[\nabla_{X},D^+]=-\frac{1}{2}\mathrm{Ric}(X^+)\cdot+\frac{1}{2}\big(X^+\lrcorner F^{1,1}_A\big)\cdot +\frac{1}{2}\big(X^-\lrcorner F_A^{0,2}\big)\cdot,
\\
\label{nablax-}
[\nabla_{X},D^-]=-\frac{1}{2}\mathrm{Ric}(X^-)\cdot+\frac{1}{2}\big(X^-\lrcorner F^{1,1}_A\big)\cdot+\frac{1}{2}\big(X^+\lrcorner F_A^{2,0}\big)\cdot,
\end{gather}
where the $2$-form $F_A$ is decomposed as $F_A=F_A^{2,0}+F_A^{1,1}+F_A^{0,2}$, into forms of type $(2,0)$, $(1,1)$, respectively $(0,2)$.
The identity \eqref{nablax} is obtained from the following straightforward computation
\begin{gather*}
\nabla_X(D\varphi) =\nabla_X( e_j\cdot \nabla_{e_j}\varphi)= e_j\cdot \mathcal{R}^A_{X,e_j}\varphi+ e_j\cdot \nabla_{e_j}\nabla_X\varphi\\
 \hphantom{\nabla_X(D\varphi)}{} \overset{\eqref{ricident}}{=}-\frac{1}{2}{\mathop{\rm Ric}}(X)\cdot\varphi+\frac{1}{2}(X\lrcorner F_A)\cdot\varphi+D(\nabla_X \varphi).
\end{gather*}

The identity \eqref{nablax+} follows from the identities
\begin{gather*}
\nabla_{X^+}(D^+\varphi) =\nabla_{X^+}( e_i^+\cdot \nabla_{e_i^-}\varphi)= e_i^+\cdot \mathcal{R}^A_{X^+,e_i^-}\varphi+ e_i^+\cdot \nabla_{e_i^-}\nabla_{X^+}\varphi\\
\hphantom{\nabla_{X^+}(D^+\varphi)}{}
=-\frac{1}{2}\mathrm{Ric}(X^+)\cdot\varphi+\frac{1}{2}\big(X^+\lrcorner F^{1,1}_A\big)\cdot\varphi+D^+(\nabla_{X^+} \varphi),
\\
\nabla_{X^-}(D^+\varphi) =\nabla_{X^-}\big( e_i^+\cdot \nabla_{e_i^-}\varphi\big)= e_i^+\cdot R_{X^-,e_i^-}\varphi+ e_i^+\cdot \nabla_{e_i^-}\nabla_{X^-}\varphi\\
\hphantom{\nabla_{X^-}(D^+\varphi)}{}
=\frac{1}{2}\big(X^-\lrcorner F^{0,2}_A\big)\cdot\varphi+D^+(\nabla_{X^-} \varphi).
\end{gather*}
The identity \eqref{nablax-} follows either by an analogous computation or by  conjugating~\eqref{nablax+}.

On a K\"ahler manifold $(M,g,J)$  endowed with any ${\mathop{\rm spin}^c}$ structure,  a spinor of the form $\varphi_r+\varphi_{r+1}\in \Gamma(\Sigma_r M\oplus \Sigma_{r+1} M)$, for some $0\leq r\leq m$, is called a {\it   K\"ahlerian Killing ${\mathop{\rm spin}^c}$ spinor} if  there exists a non-zero real constant $\alpha$, such that the following equations are satisf\/ied, for all vector f\/ields $X$,
 \begin{gather}\label{KKSSdefinition}
   \nabla_X\varphi_r=  \alpha  X^-\cdot\varphi_{r+1},\qquad
   \nabla_X\varphi_{r+1} =  \alpha    X^+\cdot\varphi_{r}.
 \end{gather}
K\"ahlerian Killing spinors lying in $\Gamma(\Sigma_{m} M\oplus \Sigma_{m+1} M) = \Gamma(\Sigma_m M)$  or in  $\Gamma(\Sigma_{-1} M\oplus \Sigma_{0} M) = \Gamma(\Sigma_0 M)$ are just parallel spinors. A direct computation shows that each  K\"ahlerian Killing ${\mathop{\rm spin}^c}$ spinor is an eigenspinor of the square of the Dirac operator. More precisely, the following equalities hold
 \begin{gather}\label{kkseigen1}
   D\varphi_r= -2(r+1)\alpha  \varphi_{r+1},\qquad D\varphi_{r+1}= -2(m-r)\alpha \varphi_{r},
 \end{gather}
which further yield
\begin{gather}\label{kkseigen2}
D^2\varphi_r=4(m-r)(r+1)\alpha^2\varphi_r,\qquad D^2\varphi_{r+1}=4(m-r)(r+1)\alpha^2\varphi_{r+1}.
 \end{gather}
In \cite{hmu}, the authors gave examples of ${\mathop{\rm spin}^c}$ structures on compact K\"ahler--Einstein  manifolds of positive scalar curvature, which carry  K\"ahlerian Killing ${\mathop{\rm spin}^c}$ spinors lying in $\Sigma_r M\oplus \Sigma_{r+1} M$, for $r\neq \frac{m\pm1}{2}$, in contrast to the spin case, where K\"ahlerian Killing spinors may only exist for $m$ odd in the middle of the decomposition \eqref{decomp}. We brief\/ly describe  these ${\mathop{\rm spin}^c}$ structures here. If the f\/irst Chern class $c_1(K_M)$ of the canonical bundle of
the K\"ahler $M$ is a non-zero cohomology class,  the greatest number $p\in \mathbb N^*$ such that
\begin{gather*}
\frac 1p c_1 (K_M) \in H^2 (M, \mathbb Z),
\end{gather*}
is called the {\it $($Fano$)$ index} of the manifold $M$. One can thus consider a $p$-th root of the canonical bundle $K_M$, i.e.,
 a complex line bundle $\mathcal L$, such that $\mathcal L^ p = K_M$. In \cite{hmu}, O.~Hijazi, S.~Montiel and F.~Urbano proved the following:

\begin{thm}[\protect{\cite[Theorem~14]{hmu}}] \label{KKSS}
Let $M$ be a $2m$-dimensional K\"ahler--Einstein compact mani\-fold with scalar
curvature $4m(m+1)$ and index $p \in  \mathbb N^*$. For each $0 \leq r \leq m+1$, there exists on $M$ a ${\mathop{\rm spin}^c}$ structure with auxiliary line bundle given by $\mathcal L^q$, where  $q = \frac{p}{m+1} (2r-m-1) \in \mathbb Z$, and carrying a K\"ahlerian Killing spinor $\psi_{r-1} + \psi_r \in \Gamma(\Sigma_{r-1} M \oplus \Sigma_{r} M)$, i.e.,  it satisfies the first-order system
\begin{gather*}
\nabla_X \psi_{r} = - X^+ \cdot \psi_{r-1},\qquad
\nabla_X \psi_{r-1}  = - X^-  \cdot \psi_{r},
\end{gather*}
for all $X \in \Gamma(TM)$.
\end{thm}

For example, if $M$ is the complex projective space $\mathbb C P^m$ of complex dimension $m$, then $p = m+1$ and $\mathcal L$ is just the tautological line bundle. We f\/ix $0 \leq r \leq m+1$ and we endow~$\mathbb C P^m$ with the ${\mathop{\rm spin}^c}$ structure whose auxiliary line bundle  is given by $\mathcal L^q$ where $q = \frac{p}{m+1} (2r-m-1) = 2r-m-1 \in \mathbb Z$. For this ${\mathop{\rm spin}^c}$ structure, the space of  K\"ahlerian Killing spinors   in $\Gamma(\Sigma_{r-1}M \oplus \Sigma_rM)$ has dimension $\binom{m+1}{r}$. A K\"ahler manifold carrying a complex contact structure necessarily has odd complex dimension $m = 2\ell+1$ and its index $p$ equals $\ell+1$.  We f\/ix $0 \leq r \leq m+1$ and we endow $M$ with the ${\mathop{\rm spin}^c}$ structure whose auxiliary line bundle is  given by $\mathcal L^q$ where $q = \frac{p}{m+1} (2r-m-1) = r-\ell-1 \in \mathbb Z$. For this ${\mathop{\rm spin}^c}$ structure, the space of  K\"ahlerian Killing spinors   in $\Gamma(\Sigma_{r-1} M\oplus \Sigma_{r}M)$ has dimension $1$.  In these examples, for $r=0$ (resp. $r = m+1$), we get the canonical (resp. anticanonical) ${\mathop{\rm spin}^c}$ structure for which K\"ahlerian Killing spinors are just parallel spinors.

\section[Eigenvalue estimates for the ${\mathop{\rm spin}^c}$ Dirac operator on K\"ahler--Einstein manifolds]{Eigenvalue estimates for the $\boldsymbol{{\mathop{\rm spin}^c}}$ Dirac operator\\ on K\"ahler--Einstein manifolds}\label{eigen}

In this section, we give a lower bound for the eigenvalues of the ${\mathop{\rm spin}^c}$ Dirac operator on a K\"ahler--Einstein manifold endowed with particular ${\mathop{\rm spin}^c}$ structures. More precisely, let $(M^{2m},g,J)$ be a~compact K\"ahler--Einstein manifold of index $p\in\mathbb{N}^*$ and of positive scalar curvature $S$, endowed with the ${\mathop{\rm spin}^c}$ structure given by $\mathcal{L}^q$, where $\mathcal L$ is the $p$-$th$ root of the canonical bundle and $q+p\in2\mathbb{Z}$ (among all powers $\mathcal L^q$,
only those satisfying $p+q \in 2\mathbb Z$ provide us a  ${\mathop{\rm spin}^c}$ structure,  cf.~\cite[Section~7]{hmu}). The curvature form $F_A$ of the induced connection $A$ on $\mathcal{L}^q$  acts on the spinor bundle as $\frac{q}{p}i\rho$. Since $(M^{2m},g,J)$ is K\"ahler--Einstein, it follows that $\rho=\frac{S}{2m}\Omega$, where $\Omega$ is the K\"ahler form. Hence,
for each $0\leq r\leq m$, we have
\begin{gather}\label{fa-ke}
(S +2 F_A )\cdot \varphi_r =\left(1-\frac{q}{p}\cdot\frac{2r-m}{m}\right) S \varphi_r, \qquad \forall\, \varphi_r\in\Gamma(\Sigma_{r}M).
\end{gather}
Let us denote by $c_r:=1-\frac{q}{p}\cdot \frac{2r-m}{m}$ and
\begin{gather*}
  a_1\colon \ \{0,\dots,m\} \to \mathbb{R}, \qquad   a_1(r):=\frac{r+1}{2r+1} c_r,\\ 
  a_2\colon \ \{0,\dots,m\} \to \mathbb{R},\qquad    a_2(r):=\frac{m-r+1}{2m-2r+1}c_r.
\end{gather*}

With the above notation, the following result holds:
\begin{prop}\label{estimgen}
Each eigenvalue $\lambda_r$ of $D^2$ restricted to $\Sigma_{r}M$ with associated eigenspinor~$\varphi_r$ satisfies the inequality
\begin{gather}\label{ineggen}
\lambda_r\geq \max \big (\min \big (a_1(r), a_1(r-1)\big ), \min\big( a_2(r), a_2(r+1)\big)\big)\cdot \frac{S}{2}.
\end{gather}
Moreover, the equality case is characterized as follows:
\begin{enumerate}[$a)$]\itemsep=0pt
 \item $D^2\varphi_r=a_1(r)\frac{S}{2}\varphi_r \Longleftrightarrow T_r\varphi_r=0$, $D^-\varphi_r=0$;
 \item $D^2\varphi_r=a_1(r-1)\frac{S}{2}\varphi_r \Longleftrightarrow  T_{r-1}(D^-\varphi_r)=0$;
 \item $D^2\varphi_r=a_2(r)\frac{S}{2}\varphi_r \Longleftrightarrow T_r\varphi_r=0$, $D^+\varphi_r=0$;
 \item $D^2\varphi_r=a_2(r+1)\frac{S}{2}\varphi_r \Longleftrightarrow  T_{r+1}(D^+\varphi_r)=0$.
\end{enumerate}
\end{prop}

\begin{proof}
For $0\leq r\leq m$ we have: $(S +2 F_A )\cdot \varphi_r = c_r S  \varphi_r$, $\forall\varphi_r\in\Gamma(\Sigma_{r}M)$. Let $r\in\{0,\dots,m\}$ be f\/ixed, $\lambda_r$ be an eigenvalue of $D^2|_{\Sigma_r M}$ and $\varphi_r\in\Gamma(\Sigma_{r}M)$ be an eigenspinor: $D^2\varphi_r=\lambda_r\varphi_r$. We distinguish two cases.

i) If $D^-\varphi_r=0$, then $|D\varphi_r|^2=|D^+\varphi_r|^2$ and \eqref{ineg} implies
\begin{gather*}
|\nabla\varphi_r|^2\geq \frac{1}{2(r+1)}|D^+\varphi_r|^2=\frac{1}{2(r+1)}|D\varphi_r|^2.
\end{gather*}
By Lemma~\ref{alg1},  it follows that
\begin{gather*}\lambda_r\geq\frac{r+1}{2(2r+1)}c_rS. \end{gather*}

 ii) If $D^-\varphi_r\neq 0$, then we consider $\varphi_r^-:= D^-\varphi_r$, which satisf\/ies $D^2\varphi_r^-=\lambda_r\varphi_r^-$ and $D^-\varphi_r^-=0$,  so in particular $|D\varphi_r^-|^2=|D^+\varphi_r^-|^2$. We now apply the argument in i) to $\varphi_r^-\in\Gamma(\Sigma_{r-1}M)$. By~\eqref{ineg}, it follows that
\begin{gather*}
|\nabla\varphi_r^-|^2\geq \frac{1}{2r}|D^+\varphi_r^-|^2=\frac{1}{2r}|D\varphi_r^-|^2.
\end{gather*}
Applying again Lemma~\ref{alg1}, we obtain $\lambda_r\geq\frac{r}{2(2r-1)} c_{r-1} S$.

Hence, we  have showed that $\lambda_r\geq \min\big (a_1(r), a_1(r-1)\big)\frac{S}{2}$. The same argument applied to the cases when $D^+\varphi_r=0$ and $D^+\varphi_r\neq 0$ proves the inequality $\lambda_r\geq \min\big (a_2(r), a_2(r+1)\big ) \frac{S}{2}$. Altogether we obtain the estimate in Proposition~\ref{estimgen}. The characterization of the equality cases is a direct consequence of Lemma~\ref{alg1}, identity~\eqref{ident} and the description of  the limiting case of inequality~\eqref{ineg}.
\end{proof}

\begin{Remark}
The inequality \eqref{ineggen} can be expressed more explicitly, by determining the maximum according to several possible cases. However, since in the sequel we will ref\/ine this eigenvalue estimate, we are only interested in the characterization of the limiting cases, which will be used later in the proof of  the equality case of the estimate (\ref{global}).
\end{Remark}
In order to ref\/ine the estimate (\ref{ineggen}), we start by the following two lemmas.

\begin{lem} \label{lemmacit}
Let $(M^{2m},g,J)$ be a compact K\"ahler--Einstein manifold of index $p$ and of positive scalar curvature $S$, endowed with a ${\mathop{\rm spin}^c}$ structure given by $\mathcal{L}^q$, where  $q+p\in2\mathbb{Z}$. For any spinor field $\varphi$ and any vector field $X$, the ${\mathop{\rm spin}^c}$ Ricci identity is given by
\begin{gather} \label{kericspinc}
e_j \cdot\mathcal{R}^A_{e_j,X} \varphi = \frac{1}{2}\mathrm{Ric}(X) \cdot\varphi- \frac{S}{4m}\frac{q}{p} (X\lrcorner i\Omega)\cdot\varphi,
\end{gather}
and it can be refined as follows
\begin{gather}\label{kericspinc+-}
e^-_j \cdot\mathcal{R}^A_{e^+_j,X^-} \varphi =\frac{1}{2}\mathrm{Ric}(X^-)\cdot \varphi -\frac{S}{4m}\frac{q}{p} X^-\cdot\varphi,
\\
\label{kericspinc-+}
e^+_j \cdot\mathcal{R}^A_{e^-_j,X^+} \varphi  = \frac{1}{2}\mathrm{Ric}(X^+)\cdot \varphi +\frac{S}{4m}\frac{q}{p} X^+\cdot\varphi.
\end{gather}
\end{lem}

\begin{proof} Since the  curvature form $F_A$ of the ${\mathop{\rm spin}^c}$ structure acts on the spinor bundle as $\frac{q}{p}i\rho=\frac{q}{p} \frac{S}{2m} i\Omega$, \eqref{kericspinc} follows directly from the Ricci identity \eqref{ricident}. The ref\/ined identities \eqref{kericspinc+-} and \eqref{kericspinc-+}  follow by replacing $X$ in \eqref{kericspinc} with $X^-$, respectively $X^+$, which is possible since both sides of the identity are complex linear in $X$, and by taking into account that when decomposing $e_j=e_j^+ + e^-_j$, the following identities (and their analogue for $X^+$) hold: $e_j \cdot\mathcal{R}^A_{e^-_j,X^-}=0$ and $e^+_j \cdot\mathcal{R}^A_{e^+_j,X^-}=0$. These last two identities are a consequence of the $J$-invariance of the curvature tensor,  i.e.,  $\mathcal{R}^A_{JX,JY}=\mathcal{R}^A_{X,Y}$,
for all vector f\/ields~$X$,~$Y$, as this implies $\mathcal{R}^A_{e^-_j,X^-}=\mathcal{R}^A_{Je^-_j,JX^-}=(-i)^2\mathcal{R}^A_{e^-_j,X^-}$ and also $e^+_j \cdot\mathcal{R}^A_{e^+_j,X^-}=Je^+_j \cdot\mathcal{R}^A_{Je^+_j,X^-}=i^2 e^+_j \cdot\mathcal{R}^A_{e^+_j,X^-}$, so they both vanish. In order to obtain the second term on the right hand side of~\eqref{kericspinc+-} and~\eqref{kericspinc-+}, we use the following identities
 of endomorphisms of the spinor bundle: $X^-\lrcorner i\Omega=X^-$ and $X^+\lrcorner i\Omega=-X^+$.
\end{proof}

\begin{lem}\label{refinedd} Under \looseness=-1 the same assumptions as in Lemma~{\rm \ref{lemmacit}}, the  refined  Schr\"odinger--Lichne\-ro\-wicz formula for ${\mathop{\rm spin}^c}$ K\"ahler manifolds  for the action on each eigenbundle $\Sigma_{r}M$ is given by
\begin{gather}\label{sl1}
 2\nabla^{{1,0}^*}\nabla^{1,0}=D^2-\frac{S}{4}-\frac{i}{2}\rho-\frac{m-r}{2m}\frac{q}{p}S,
\\
\label{sl2}
 2\nabla^{{0,1}^*}\nabla^{0,1}=D^2-\frac{S}{4}+\frac{i}{2}\rho+\frac{r}{2m}\frac{q}{p}S,
\end{gather}
where $\nabla^{1, 0}$  $($resp.\ $\nabla^{0,1})$  is the holomorphic $($resp.\ antiholomorphic$)$  part of $\nabla$, i.e., the projections of~$\nabla$ onto the following two components
\begin{gather*}
\nabla\colon \ \Gamma(\Sigma_{r}M)\to \Gamma(\Lambda^{1,0}M\otimes\Sigma_{r}M)\oplus \Gamma(\Lambda^{0,1}M\otimes\Sigma_{r}M).
\end{gather*} They are locally defined, for all vector fields~$X$,  by
\begin{gather*}
\nabla^{1, 0}_X  = g(X, e_i^-) \nabla_{e_i^+}=\nabla_{X^+} \qquad \text{and}\qquad \nabla^{0, 1}_X = g(X, e_i^+) \nabla_{e_i^-}=\nabla_{X^-},
\end{gather*}
where $\{e_1,\dots, e_{2m}\}$ is a local orthonormal basis of~$TM$.
\end{lem}

\begin{proof}
Let $\{e_1,\dots, e_{2m}\}$ be a local orthonormal basis of $TM$ (identif\/ied with $\Lambda^1 M$ via the met\-ric~$g$), parallel at the point where the computation is made. We recall that the formal adjoints $\nabla^{{1,0}^*}$ and $\nabla^{{1,0}^*}$ are given by the following formulas (for a proof, see, e.g., \cite[Lemma~20.1]{am_lectures})
\begin{gather*}
\nabla^{{1,0}^*}\colon \ \Gamma(\Lambda^{1,0}M\otimes\Sigma_{r}M) \longrightarrow \Gamma(\Sigma_{r}M), \qquad \nabla^{{1,0}^*}(\alpha\otimes \varphi)=(\delta\alpha)\varphi-\nabla_
{\alpha}\varphi,\\
 \nabla^{{0,1}^*}\colon \ \Gamma(\Lambda^{0,1}M\otimes\Sigma_{r}M) \longrightarrow \Gamma(\Sigma_{r}M), \qquad \nabla^{{0,1}^*}(\alpha\otimes \varphi)=(\delta\alpha)\varphi-\nabla_
{\alpha}\varphi.
\end{gather*}
We thus obtain for the corresponding Laplacians
\begin{gather}\label{laplac}
\nabla^{{1,0}^*}\nabla^{1,0}\varphi
=\nabla^{{1,0}^*}(e^-_j\otimes \nabla_{e^+_j}\varphi)=- \nabla_{e^-_j}\nabla_{e^+_j},
\end{gather}
since $\delta e^-_j=0$, as the basis is  parallel at the given point, and $g(\cdot, e^-_j)\in\Lambda^{1,0}M$. Analogously, or by conjugation, we have $\nabla^{{0,1}^*}\nabla^{0,1}\varphi
=- \nabla_{e^+_j}\nabla_{e^-_j}$.
We now prove~\eqref{sl1}. By a similar computation, one obtains \eqref{sl2}
\begin{gather*}
 2\nabla^{{1,0}^*}\nabla^{1,0} \overset{\eqref{laplac}}{=}-2g(e_i,e_j)\nabla_{e_i^-}\nabla_{e_j^+}=(e_i\cdot e_j + e_j\cdot e_i) \cdot \nabla_{e_i^-}\nabla_{e_j^+}\\
\hphantom{2\nabla^{{1,0}^*}\nabla^{1,0}}{}
  =D^+D^- + e_j\cdot e_i\cdot (\nabla_{e_j^+}\nabla_{e_i^-} - R_{e_j^+,e_i^-})=D^+D^- + D^-D^+ + e_j^-\cdot e^+_i\cdot R_{e_i^-, e_j^+}\\
 \hphantom{2\nabla^{{1,0}^*}\nabla^{1,0}}{}
 \overset{\eqref{kericspinc-+}}{=}D^2+e^-_j\cdot \left( \frac{1}{2}\mathrm{Ric}(e_j^+)+\frac{S}{4m}\frac{q}{p}e_j^+\right)\\
  \hphantom{2\nabla^{{1,0}^*}\nabla^{1,0}}{}
  \overset{\eqref{kecontr1},~\eqref{kecontr3}}{=}D^2 -\frac{1}{2}\left(\frac{S}{2}+i\rho\right)-\frac{m-r}{2m}\frac{q}{p}S.\tag*{\qed}
 \end{gather*}
 \renewcommand{\qed}{}
\end{proof}

\begin{thm}\label{estimke}
Let $(M^{2m}, g, J)$ be a compact K\"ahler--Einstein manifold of index $p$ and positive scalar curvature $S$, carrying the ${\mathop{\rm spin}^c}$ structure given by $\mathcal{L}^q$ with $q+p\in 2\mathbb{Z}$, where $\mathcal{L}^p=K_M$. We assume that  $p\geq |q|$.
Then, for each $r\in\{0, \dots, m\}$, any eigenvalue $\lambda_r$ of $D^2|_{\Gamma(\Sigma_r M )}$ satisfies the inequality
\begin{gather}\label{eqestimke}
 \lambda_r\geq e(r)\frac{S}{2},
 \end{gather}
where
\begin{gather*}
e\colon \ [0,m]\to\mathbb{R}, \qquad  e(x)=\begin{cases}
   \displaystyle    e_1(x)=\frac{m-x}{m}\left(1+\frac{q}{p}\right), & \displaystyle\text{if} \ \ x\leq \left(1+\frac{q}{p}\right)\frac{m}{2},
   \vspace{1mm}\\
    \displaystyle    e_2(x)=\frac{x}{m}\left(1-\frac{q}{p}\right),  & \displaystyle\text{if} \ \ x\geq \left(1+\frac{q}{p}\right)\frac{m}{2}.                                                                                                                        \end{cases}
\end{gather*}
Moreover, equality is attained if and only if the corresponding eigenspinor $\varphi_r\in\Gamma(\Sigma_r M)$ is an antiholomorphic spinor: $\nabla^{1,0}\varphi_r=0$, if $r\leq \left(1+\frac{q}{p}\right)\frac{m}{2}$, respectively a holomorphic spinor: $\nabla^{0,1}\varphi_r=0$, if $r \geq \big(1+\frac{q}{p}\big)\frac{m}{2}$.
\end{thm}

\begin{proof}
First we notice that our assumption $|q|\leq p$ implies that the lower bound in~\eqref{eqestimke} is non-negative and that $0\leq \big(1+\frac{q}{p}\big)\frac{m}{2}\leq m$. The formulas~\eqref{sl1} and~\eqref{sl2} applied to $\varphi_r$ yield, after taking the scalar product with $\varphi_r$ and integrating over $M$, the following inequalities
 \begin{gather*}\lambda_r\geq \frac{m-r}{m}\left(1+\frac{q}{p}\right)\frac{S}{2},\qquad
 \lambda_r\geq \frac{r}{m}\left(1-\frac{q}{p}\right)\frac{S}{2},\end{gather*}
and equality is attained if and only if the corresponding eigenspinor $\varphi_r$ satisf\/ies $\nabla^{1,0}\varphi_r=0$, resp.~$\nabla^{0,1}\varphi_r=0$. Hence, for any $0\leq r\leq m$ we obtain the following lower bound:
\begin{gather*}
\lambda_r\geq \max\left(\frac{m-r}{m}\left(1+\frac{q}{p}\right),\frac{r}{m}\left(1-\frac{q}{p}\right)\right)=  e(r)\frac{S}{2}.\tag*{\qed}
\end{gather*}
\renewcommand{\qed}{}
\end{proof}

\begin{Remark}\label{comp} Let us denote $\frac{q}{p}\cdot\frac{m+1}{2}+\frac{m-1}{2}$  by $b$.
Comparing the estimate given by Theorem~\ref{estimke} with the estimate from Proposition~\ref{estimgen}, we obtain for $r\leq  b$
\begin{gather*}e(r)-a_1(r)=\frac{(m+1)\frac{q}{p}+m-1-2r}{m(2r+1)}=-\frac{2(r-b)}{m(2r+1)}.\end{gather*}
Hence, for $r \leq  b$, we have $e(r)\geq a_1(r)$ and $e(r)=a_1(r)$ if\/f $r=b\in\mathbb{N}$.
Similarly, for $r\geq b+1$, we compute
\begin{gather*}e(r)-a_2(r)=\frac{2(m-r)(r-b-1)}{m(2m-2r+1)}.\end{gather*}
Hence, for  $r\geq b+1$, we have $e(r)\geq a_2(r)$ and $e(r)=a_2(r)$ if\/f $r=b+1\in\mathbb{N}$.
\end{Remark}

Theorem~\ref{estimke} implies the global lower bound for the eigenvalues of the ${\mathop{\rm spin}^c}$ Dirac operator acting on the whole spinor bundle in Theorem \ref{globestimke}. We are now ready to prove this result.

\begin{proof}[Proof of Theorem~\ref{globestimke}]
\hspace{0.01cm} Since the lower bound established in Theorem~\ref{estimke} decreases on $\big(0,\big(1+\frac{q}{p}\big)\frac{m}{2}\big)$ and increases on $\big(\big(1+\frac{q}{p}\big)\frac{m}{2},m\big)$, we obtain the following global estimate
\begin{gather*}\lambda\geq e\left(\left(1+\frac{q}{p}\right)\frac{m}{2}\right)=\frac{1}{2}\left(1-\frac{q^2}{p^2}\right)\frac{S}{2}.\end{gather*}

However, this estimate is not sharp. Otherwise, this would imply that $\big(1+\frac{q}{p}\big)\frac{m}{2}\in\mathbb{N}$ and the limiting eigenspinor would be, according to the characterization of the equality case in Theo\-rem~\ref{estimke}, both holomorphic and antiholomorphic, hence parallel and, in particular, harmonic. This fact together with the Lichnerowicz--Schr\"odinger formula~\eqref{sl} and the fact that the scalar curvature is positive leads to a contradiction.

We now assume that there exists an $r\in\mathbb{N}$, such that $b<r< \big(1+\frac{q}{p}\big)\frac{m}{2}$ and the equality in~\eqref{eqestimke} is attained. We obtain a contradiction as follows. Let $\varphi_r$ be the corresponding eigenspinor: $D^2\varphi_r=e_1(r)\frac{S}{2}\varphi_r$ and $\nabla^{1,0}\varphi_r=0$. Then $D^+\varphi_r\in\Sigma_{r+1}M$ is also an eigenspinor of $D^2$ to the eigenvalue $e_1(r)\frac{S}{2}$ (note that $D^+\varphi_r\neq 0$, otherwise $\varphi_r$ would be a harmonic spinor and we could conclude as above). However, for all $r>b$, the strict inequality $e_2(r+1)>e_1(r)$ holds. Since $r+1>\big(1+\frac{q}{p}\big)\frac{m}{2}$, this contradicts the estimate \eqref{eqestimke}. The same argument as above shows that there exists no $r\in\mathbb{N}$, such that $\big(1+\frac{q}{p}\big)\frac{m}{2}<r<b+1 $ and the equality in \eqref{eqestimke} is attained. Hence, we obtain the following global estimate
\begin{gather*}
\lambda\geq e_1(b)\frac{S}{2}=e_2(b+1)\frac{S}{2}=\frac{m+1}{2m}\left(1-\frac{q^2}{p^2}\right)\frac{S}{2}=\left(1-\frac{q^2}{p^2}\right)(m+1)^2.
\end{gather*}
According to Theorem~\ref{estimke}, the equality is attained if and only if $b \in \mathbb N$ and the corresponding eigenspinors  $\varphi_{b}\in\Gamma(\Sigma_{b} M)$ and $\varphi_{b+1}\in\Gamma(\Sigma_{b+1}M)$ to the eigenvalue $\big(1-\frac{q^2}{p^2}\big)(m+1)^2$ are antiholomorphic resp. holomorphic spinors: $\nabla^{1,0}\varphi_{b}=0$, $\nabla^{0,1}\varphi_{b+1}=0$. In particular, this implies $D^-\varphi_{b}=0$ and $D^+\varphi_{b+1}=0$. By Remark~\ref{comp}, we have: $e_1(b)=a_1(b)$ and $e_2(b+1)=a_2(b+1)$. Hence, the characterization of the equality case in Proposition~\ref{estimgen} yields $T_{b}\varphi_{b}=0$ and $T_{b+1}\varphi_{b+1}=0$, which further imply
\begin{gather}\label{twis1}
\nabla_X \varphi_{b}=-\frac{1}{2(b+1)}X^-\cdot D^+\varphi_{b}=-\frac{1}{2(b+1)}X^-\cdot D\varphi_{b},
\\
\nabla_X \varphi_{b+1}=-\frac{1}{2(m-b)}X^+\cdot D^-\varphi_{b+1}=-\frac{1}{2(m-b)}X^+\cdot D\varphi_{b+1}.\nonumber
\end{gather}

We now show that the spinors $\varphi_{b}+\frac{1}{(m+1)\left(1+\frac{q}{p}\right)}D\varphi_{b}\in\Gamma(\Sigma_{b} M\oplus \Sigma_{b+1} M)$ and $\varphi_{b+1}+\frac{1}{(m+1)\left(1-\frac{q}{p}\right)}D\varphi_{b+1}\in\Gamma(\Sigma_{b+1} M\oplus \Sigma_{b} M)$ are K\"ahlerian Killing  ${\mathop{\rm spin}^c}$ spinors. Note that for $q=0$ (corresponding to the spin case), it follows that $\varphi_{b}+\frac{1}{m+1}D\varphi_{b}, \varphi_{b+1}+\frac{1}{m+1}D\varphi_{b+1}\in\Gamma(\Sigma_{b} M\oplus \Sigma_{b+1} M)$ are eigenspinors of the Dirac operator corresponding to the smallest possible eigenvalue $m+1$, i.e., K\"ahlerian Killing spinors. From \eqref{twis1} it follows
\begin{gather}\label{twiss1}
\nabla_X \varphi_{b}=-X^-\cdot \frac{1}{(m+1)\left(1+\frac{q}{p}\right)} D\varphi_{b}.
\end{gather}
Applying \eqref{nablax+}  to $\varphi_{b}$ in this case for ${\mathop{\rm Ric}}=\frac{S}{2m}g=2(m+1)g$ and $F_A=\frac{q}{p}\frac{S}{2m}i\Omega=2(m+1)\frac{q}{p}i\Omega$, we get
\begin{gather}\label{twiss2}
\nabla_X (D^+\varphi_{b})=-(m+1)\left(1+\frac{q}{p}\right)X^+\cdot \varphi_{b}.
\end{gather}
 According to the def\/ining  equation \eqref{KKSSdefinition} of a K\"ahlerian Killing ${\mathop{\rm spin}^c}$ spinor, equa\-tions~\eqref{twiss1} and~\eqref{twiss2} imply that the spinor $\varphi_{b}+\frac{1}{(m+1)\big(1+\frac{q}{p}\big)}D\varphi_{b}\in\Gamma(\Sigma_{b} M\oplus \Sigma_{b+1} M)$ is a  K\"ahlerian Killing~${\mathop{\rm spin}^c}$ spinor. A similar computation yields that $\varphi_{b+1}+\frac{1}{(m+1)\big(1-\frac{q}{p}\big)}D\varphi_{b+1}$ is a K\"ahlerian Killing~${\mathop{\rm spin}^c}$  spinor.  Conversely, if $\varphi_{b}+\varphi_{b+1}\in \Gamma(\Sigma_{b}M\oplus \Sigma_{b+1} M)$ is a K\"ahlerian Killing~${\mathop{\rm spin}^c}$ spinor, then according to~\eqref{kkseigen2}, $\varphi_{b}$ and $\varphi_{b+1}$ are eigenspinors of $D^2$ to the eigenvalue $4(m-b)(b+1)=\big(1-\frac{q^2}{p^2}\big)(m+1)^2$. This concludes the proof.
\end{proof}

\begin{Remark}
If $q=0$, which corresponds to the spin case, the assumption $p \geq \vert q \vert = 0$ is trivial and we recover from Theorem~\ref{estimke} and Theorem~\ref{globestimke} Kirchberg's estimates on K\"ahler--Einstein spin manifolds: the lower bound \eqref{kirchoddeven} for $m$ odd, namely $\lambda^2\geq \frac{m+1}{4m}S=e\big(\frac{m+1}{2}\big)\frac{S}{2}$,  and the lower bound~\eqref{kirchke} for $m$ even, namely $\lambda^2\geq \frac{m+2}{4m}S=e\big(\frac{m}{2}+1\big)\frac{S}{2}$. In the latter case, when~$m$ is even, the equality in \eqref{global} cannot be attained, as $b=\frac{m}{2}-\frac{1}{2}\notin\mathbb{N}$. Also for $r=\frac{m}{2}$ the inequality~\eqref{eqestimke} is strict, since otherwise it would imply, according to the characterization of the equality case in Theorem~\ref{estimke}, that the corresponding eigenspinor $\varphi\in\Sigma_{\frac{m}{2}}M$ is parallel, in contradiction to the positivity of the scalar curvature. Note that the same argument as in the proof of Theorem~\ref{globestimke} shows that there cannot exist an eigenspinor~$\varphi\in \Sigma_{\frac{m}{2}}M$ of $D^2$ to an eigenvalue strictly smaller than the lowest bound for $r=\frac{m}{2}\pm 1$, since otherwise~$D^+\varphi$ and~$D^-\varphi$ would either be eigenspinors or would vanish, leading in both cases to a contradiction.
Hence, from the estimate~\eqref{eqestimke} and the fact that the function~$e_1$ decreases on $\big(0,\big(1+\frac{q}{p}\big)\frac{m}{2}\big)$
and~$e_2$ increases on $\big(\big(1+\frac{q}{p}\big)\frac{m}{2},m\big)$, it follows that the lowest possible bound for $\lambda^2$ in this case is given by $e_1\big(\frac{m}{2}-1\big)S=e_2\big(\frac{m}{2}+1\big)S=\frac{m+2}{4m}$. If $q=- p$ (resp.~$q = p$), which corresponds to the canonical (resp.\ anti-canonical) ${\mathop{\rm spin}^c}$ structure, the lower bound in Theorem~\ref{globestimke} equals $0$ and is attained by the parallel spinors in~$\Sigma_{0}M$ (resp.~$\Sigma_{m}M$),  cf.~\cite{Moro1}.
\end{Remark}

\section{Harmonic forms on limiting K\"ahler--Einstein manifolds}

In this section we give an application for the eigenvalue estimate of the ${\mathop{\rm spin}^c}$ Dirac operator established in Theorem~\ref{globestimke}.
Namely, we extend to ${\mathop{\rm spin}^c}$ spinors the result of A.~Moroianu \cite{am} stating that
the Clif\/ford multiplication between a harmonic  ef\/fective form of nonzero degree and a~K\"ahlerian Killing spinor vanishes.
As above, $(M^{2m}, g)$ denotes a $2m$-dimensional K\"ahler--Einstein compact manifold of  index $p$ and normalized scalar curvature $4m(m+1)$, which carries the ${\mathop{\rm spin}^c}$ structure given by $\mathcal{L}^q$ with $q+p\in 2\mathbb{Z}$, where $\mathcal{L}^p=K_M$. We call $M$ a {\it limiting mani\-fold} if equality in \eqref{global} is achieved on $M$,  which is by Theorem~\ref{globestimke} equivalent to the existence of a~K\"ahlerian Killing ${\mathop{\rm spin}^c}$ spinor in $\Sigma_{r}M\oplus \Sigma_{r+1} M$ for $r=\frac{q}{p}\cdot\frac{m+1}{2}+\frac{m-1}{2}\in\mathbb{N}$.
Let $\psi=\psi_{r-1}+\psi_{r}\in \Gamma(\Sigma_{r-1} M \oplus \Sigma_{r} M)$ be such a spinor, i.e., $\Omega \cdot \psi_{r-1} = i (2r-2-m) \psi_{r-1}$, $\Omega \cdot \psi_{r} = i (2r-m) \psi_{r}$ and the following equations are satisf\/ied
\begin{gather*}
\nabla_{X^+}\psi_{r} =- X^+ \cdot \psi_{r-1},\qquad
\nabla_{X^-} \psi_{r-1} =- X^-\cdot \psi_{r}.
\end{gather*}

By \eqref{kkseigen1},  we  have
\begin{gather*}D\psi_{r}=2(m-r+1)\psi_{r-1}, \qquad D\psi_{r-1}=2r\psi_{r}.\end{gather*}
Recall that a form $\omega$ on a K\"ahler manifold is called {\it effective} if $\Lambda \omega = 0$, where $\Lambda$ is the adjoint of the operator
$L\colon\Lambda^{*} M \longrightarrow \Lambda^{*+2} M $, $L(\omega):=\omega\wedge\Omega$. More precisely, $\Lambda$ is given by the formula:
$\Lambda = -2 \overset{2m}{\underset{j=1}{\sum}} e_j^+  \lrcorner e_j^- \lrcorner$. Moreover, one can check that
\begin{gather*}
(\Lambda L - L \Lambda ) \omega = (m-t) \omega, \qquad \forall\,  \omega \in \Lambda^t M.
\end{gather*}

\begin{lem}
Let $\psi=\psi_{r-1}+\psi_{r}\in\Gamma(\Sigma_{r-1} M \oplus \Sigma_{r} M)$ be a  K\"ahlerian Killing  ${\mathop{\rm spin}^c}$ spinor and~$\omega$ a harmonic effective form of type $(k, k')$. Then, we have
\begin{gather}\label{domeg1}
D(\omega\cdot\psi_{r})=2(-1)^{k+k'} (m-r+1-k') \omega\cdot\psi_{r-1},
\\ \label{domeg2}
D(\omega\cdot\psi_{r-1})=2(-1)^{k+k'} (r-k) \omega\cdot\psi_{r}.
\end{gather}
\end{lem}

\begin{proof}
The following general formula  holds for {\color{black}{any form}} $\omega$ of degree $\text{deg}(\omega)$ and any spinor $\varphi$
\begin{gather*}D(\omega\cdot\varphi)=(d\omega+\delta\omega)\cdot\varphi+(-1)^{\text{deg}(\omega)} \omega\cdot D\varphi -2\sum_{j=1}^{2m}(e_j\lrcorner\omega)\cdot \nabla_{e_j}\varphi.\end{gather*}
Applying this formula to an ef\/fective harmonic form $\omega$ of type $(k, k')$ and to the components of the  K\"ahlerian Killing ${\mathop{\rm spin}^c}$ spinor $\psi$, we obtain
\begin{gather*}
D(\omega\cdot\psi_{r}) = (-1)^{k+k'} \omega\cdot D\psi_{r-1}-2\sum_{j=1}^{2m}(e_j\lrcorner\omega)\cdot \nabla_{e_j}\psi_{r}\\
\hphantom{D(\omega\cdot\psi_{r})}{}
 =(-1)^{k+k'} 2(m-r+1) \omega\cdot\psi_{r-1}+2\sum_{j=1}^{2m}(e^-_j\lrcorner\omega)\cdot e^+_j\cdot \psi_{r-1}\\
\hphantom{D(\omega\cdot\psi_{r})}{}
=2(-1)^{k+k'} \left[(m-r+1) \omega\cdot\psi_{r-1}+\left(\sum_{j=1}^{2m} e^+_j\wedge(e^-_j\lrcorner\omega)\right)\cdot\psi_{r-1}\right].
\end{gather*}

Since $\omega$ is ef\/fective, we have for any spinor $\varphi$ that
\begin{gather*}(e^-_j\lrcorner\omega)\cdot e^+_j\cdot \varphi=(-1)^{k+k'-1} \big(e^+_j\wedge(e^-_j\lrcorner\omega)+e^+_j\lrcorner e^-_j\lrcorner\omega \big)\cdot \varphi.\end{gather*}
Thus, we conclude $D(\omega\cdot\psi_{r}) = 2 (-1)^{k+k'} (m-r+1-k') \omega\cdot \psi_{r-1}$.
Analogously we obtain $D(\omega\cdot\psi_{r-1})=2(-1)^{k+k'} (r-k) \omega\cdot\psi_{r}.$
\end{proof}

Now, we are able to state the main result of this section, which extends the result of A.~Moroianu mentioned in the introduction to the~${\mathop{\rm spin}^c}$ setting:
\begin{thm} \label{eff}
On a compact K\"ahler--Einstein limiting manifold, the Clifford multiplication of a~harmonic effective form of nonzero degree with the corresponding  K\"ahlerian Killing ${\mathop{\rm spin}^c}$ spinor vanishes.
\end{thm}

\begin{proof}
Equations \eqref{domeg1} and \eqref{domeg2} imply that
\begin{gather*}
D^2 (\omega \cdot \psi) = 4 (r-k)(m-r+1-k')  \omega\cdot \psi.
\end{gather*}
Note that for all values of $k,k' \in\{0,\dots,m\}$ and $r\in\{0,\dots, m+1\}$, either $4 (r-k)(m-r+1-k') \leq 0$, or $4 (r-k)(m-r+1-k') < 4r(m-r+1)$, which for $r=b+1$ is exactly the lower bound obtained in Theorem~\ref{globestimke} for the eigenvalues of~$D^2$. This shows that $\omega\cdot\psi=0$.
\end{proof}

K\"ahler--Einstein manifolds carrying a complex contact structure are examples of odd-di\-men\-sional  K\"ahler manifolds with K\"ahlerian Killing ${\mathop{\rm spin}^c}$ spinors in $\Sigma_{r-1} M \oplus\Sigma_r M $ for the ${\mathop{\rm spin}^c}$ structure (described in the introduction) whose auxiliary line bundle is given by $\mathcal L^q$ and $q = r-\ell-1$, where $m=2\ell+1$. Thus, the result of A.~Moroianu is obtained as a special case of Theorem~\ref{eff}.

\subsection*{Acknowledgments}

The f\/irst named author gratefully acknowledges the f\/inancial support of the Berlin Mathematical
School (BMS) and would like to thank the University of Potsdam, especially
Christian B\"ar and his group, for their generous support and friendly welcome during summer 2013
and summer 2014. The f\/irst named author thanks also the Faculty of Mathematics of the
University of Regensburg  for its support and hospitality during his two visits in July 2013 and July 2014.
The authors are very much indebted to Oussama Hijazi and Andrei Moroianu for many useful discussions. Both authors thank the editor and the referees for carefully reading the paper and for providing constructive comments, which substantially improved it.

\pdfbookmark[1]{References}{ref}
\LastPageEnding


\begin{thebibliography}{99}
\footnotesize \itemsep=0pt

\bibitem{AS}
Atiyah M.F., Singer I.M., The index of elliptic operators.~{I}, \href{http://dx.doi.org/10.2307/1970715}{\textit{Ann. of
  Math.}} \textbf{87} (1968), 484--530.

\bibitem{baer}
B{\"a}r C., Real {K}illing spinors and holonomy, \href{http://dx.doi.org/10.1007/BF02102106}{\textit{Comm. Math. Phys.}}
  \textbf{154} (1993), 509--521.

\bibitem{Ba}
B{\"a}r C., Extrinsic bounds for eigenvalues of the {D}irac operator,
  \href{http://dx.doi.org/10.1023/A:1006550532236}{\textit{Ann. Global Anal. Geom.}} \textbf{16} (1998), 573--596.

\bibitem{bookspin}
Bourguignon J.-P., Hijazi O., Milhorat J.-L., Moroianu A., Moroianu S.,
  A~spinorial approach to Riemannian and conformal geometry, \textit{EMS Monographs in
  Mathematics}, European Mathematical Society, 2015.

\bibitem{don}
Donaldson S.K., The {S}eiberg--{W}itten equations and {$4$}-manifold topology,
  \href{http://dx.doi.org/10.1090/S0273-0979-96-00625-8}{\textit{Bull. Amer. Math. Soc.}} \textbf{33} (1996), 45--70.

\bibitem{fr1}
Friedrich Th., Der erste {E}igenwert des {D}irac-{O}perators einer kompakten,
  {R}iemannschen {M}annigfaltigkeit nichtnegativer {S}kalarkr\"ummung,
  \href{http://dx.doi.org/10.1002/mana.19800970111}{\textit{Math. Nachr.}} \textbf{97} (1980), 117--146.

\bibitem{fr_book}
Friedrich Th., Dirac operators in {R}iemannian geometry, \href{http://dx.doi.org/10.1090/gsm/025}{\textit{Graduate
  Studies in Mathematics}}, Vol.~25, Amer. Math. Soc., Providence, RI, 2000.

\bibitem{gursky}
Gursky M.J., LeBrun C., Yamabe invariants and {${\rm Spin}^c$} structures,
  \href{http://dx.doi.org/10.1007/s000390050120}{\textit{Geom. Funct. Anal.}} \textbf{8} (1998), 965--977,
  \href{http://arxiv.org/abs/dg-ga/9708002}{dg-ga/9708002}.

\bibitem{her}
Herzlich M., Moroianu A., Generalized {K}illing spinors and conformal
  eigenvalue estimates for {${\rm Spin}^c$} manifolds, \href{http://dx.doi.org/10.1023/A:1006546915261}{\textit{Ann. Global
  Anal. Geom.}} \textbf{17} (1999), 341--370.

\bibitem{hijconf}
Hijazi O., A conformal lower bound for the smallest eigenvalue of the {D}irac
  operator and {K}illing spinors, \href{http://dx.doi.org/10.1007/BF01210797}{\textit{Comm. Math. Phys.}} \textbf{104}
  (1986), 151--162.

\bibitem{hij}
Hijazi O., Eigenvalues of the {D}irac operator on compact {K}\"ahler manifolds,
  \href{http://dx.doi.org/10.1007/BF02173430}{\textit{Comm. Math. Phys.}} \textbf{160} (1994), 563--579.

\bibitem{hmu}
Hijazi O., Montiel S., Urbano F., {${\rm Spin}^c$} geometry of {K}\"ahler
  manifolds and the {H}odge {L}aplacian on minimal {L}agrangian submanifolds,
  \href{http://dx.doi.org/10.1007/s00209-006-0936-8}{\textit{Math.~Z.}} \textbf{253} (2006), 821--853.

\bibitem{HMZ1}
Hijazi O., Montiel S., Zhang X., Dirac operator on embedded hypersurfaces,
  \href{http://dx.doi.org/10.4310/MRL.2001.v8.n2.a8}{\textit{Math. Res. Lett.}} \textbf{8} (2001), 195--208,
  \href{http://arxiv.org/abs/math.DG/0012262}{math.DG/0012262}.

\bibitem{HMZ2}
Hijazi O., Montiel S., Zhang X., Eigenvalues of the {D}irac operator on
  manifolds with boundary, \href{http://dx.doi.org/10.1007/s002200100475}{\textit{Comm. Math. Phys.}} \textbf{221} (2001),
  255--265, \href{http://arxiv.org/abs/math.DG/0012261}{math.DG/0012261}.

\bibitem{HMZ02}
Hijazi O., Montiel S., Zhang X., Conformal lower bounds for the {D}irac
  operator of embedded hypersurfaces, \textit{Asian~J. Math.} \textbf{6}
  (2002), 23--36.

\bibitem{HZ1}
Hijazi O., Zhang X., Lower bounds for the eigenvalues of the {D}irac operator.
  {I}.~{T}he hypersurface {D}irac operator, \href{http://dx.doi.org/10.1023/A:1010749808691}{\textit{Ann. Global Anal. Geom.}}
  \textbf{19} (2001), 355--376.

\bibitem{HZ2}
Hijazi O., Zhang X., Lower bounds for the eigenvalues of the {D}irac operator.
  {II}.~{T}he submanifold {D}irac operator, \href{http://dx.doi.org/10.1023/A:1011663603699}{\textit{Ann. Global Anal. Geom.}}
  \textbf{20} (2001), 163--181.

\bibitem{hit}
Hitchin N., Harmonic spinors, \href{http://dx.doi.org/10.1016/0001-8708(74)90021-8}{\textit{Adv. Math.}} \textbf{14} (1974), 1--55.

\bibitem{kirch86}
Kirchberg K.-D., An estimation for the f\/irst eigenvalue of the {D}irac operator
  on closed {K}\"ahler manifolds of positive scalar curvature, \href{http://dx.doi.org/10.1007/BF00128050}{\textit{Ann.
  Global Anal. Geom.}} \textbf{4} (1986), 291--325.

\bibitem{kirch2}
Kirchberg K.-D., The f\/irst eigenvalue of the {D}irac operator on {K}\"ahler
  manifolds, \href{http://dx.doi.org/10.1016/0393-0440(90)90001-J}{\textit{J.~Geom. Phys.}} \textbf{7} (1990), 449--468.

\bibitem{kirch}
Kirchberg K.-D., Killing spinors on {K}\"ahler manifolds, \href{http://dx.doi.org/10.1007/BF00773453}{\textit{Ann. Global
  Anal. Geom.}} \textbf{11} (1993), 141--164.

\bibitem{kircheven}
Kirchberg K.-D., Eigenvalue estimates for the {D}irac operator on
  {K}\"ahler--{E}instein manifolds of even complex dimension, \href{http://dx.doi.org/10.1007/s10455-010-9212-6}{\textit{Ann.
  Global Anal. Geom.}} \textbf{38} (2010), 273--284, \href{http://arxiv.org/abs/0912.1451}{arXiv:0912.1451}.

\bibitem{ks}
Kirchberg K.-D., Semmelmann U., Complex contact structures and the f\/irst
  eigenvalue of the {D}irac operator on {K}\"ahler manifolds, \href{http://dx.doi.org/10.1007/BF01895834}{\textit{Geom.
  Funct. Anal.}} \textbf{5} (1995), 604--618.

\bibitem{spin}
Lawson H.B., Michelsohn M.-L., Spin geometry, \textit{Princeton Mathematical
  Series}, Vol.~38, Princeton University Press, Princeton, NJ, 1989.

\bibitem{lebrun1}
LeBrun C., Einstein metrics and {M}ostow rigidity, \href{http://dx.doi.org/10.4310/MRL.1995.v2.n1.a1}{\textit{Math. Res. Lett.}}
  \textbf{2} (1995), 1--8, \href{http://arxiv.org/abs/dg-ga/9411005}{dg-ga/9411005}.

\bibitem{lebrun2}
LeBrun C., Four-manifolds without {E}instein metrics, \href{http://dx.doi.org/10.4310/MRL.1996.v3.n2.a1}{\textit{Math. Res. Lett.}}
  \textbf{3} (1996), 133--147, \href{http://arxiv.org/abs/dg-ga/9511015}{dg-ga/9511015}.

\bibitem{lich0}
Lichnerowicz A., Spineurs harmoniques, \textit{C.~R.~Acad. Sci. Paris}
  \textbf{257} (1963), 7--9.

\bibitem{am_odd}
Moroianu A., La premi\`ere valeur propre de l'op\'erateur de {D}irac sur les
  vari\'et\'es k\"ahl\'eriennes compactes, \href{http://dx.doi.org/10.1007/BF02099477}{\textit{Comm. Math. Phys.}}
  \textbf{169} (1995), 373--384.

\bibitem{am}
Moroianu A., Formes harmoniques en pr\'esence de spineurs de {K}illing
  k\"ahl\'eriens, \textit{C.~R.~Acad. Sci. Paris S\'er.~I Math.} \textbf{322}
  (1996), 679--684.

\bibitem{am_even}
Moroianu A., On {K}irchberg's inequality for compact {K}\"ahler manifolds of
  even complex dimension, \href{http://dx.doi.org/10.1023/A:1006543304443}{\textit{Ann. Global Anal. Geom.}} \textbf{15} (1997),
  235--242.

\bibitem{Moro1}
Moroianu A., Parallel and {K}illing spinors on {${\rm Spin}^c$} manifolds,
  \href{http://dx.doi.org/10.1007/s002200050142}{\textit{Comm. Math. Phys.}} \textbf{187} (1997), 417--427.

\bibitem{am_lectures}
Moroianu A., Lectures on {K}\"ahler geometry, \href{http://dx.doi.org/10.1017/CBO9780511618666}{\textit{London Mathematical
  Society Student Texts}}, Vol.~69, Cambridge University Press, Cambridge, 2007.

\bibitem{ms}
Moroianu A., Semmelmann U., K\"ahlerian {K}illing spinors, complex contact
  structures and twistor spaces, \textit{C.~R.~Acad. Sci. Paris S\'er.~I Math.}
  \textbf{323} (1996), 57--61, \href{http://arxiv.org/abs/dg-ga/9502003}{dg-ga/9502003}.

\bibitem{nakadthesis}
Nakad R., Special submanifolds of {${\rm Spin}^c$} manifolds, Ph.D.~Thesis,
  Institut \'Elie Cartan, France, 2011.

\bibitem{JRRN}
Nakad R., Roth J., Hypersurfaces of {${\rm Spin}^c$} manifolds and {L}awson
  type correspondence, \href{http://dx.doi.org/10.1007/s10455-012-9321-5}{\textit{Ann. Global Anal. Geom.}} \textbf{42} (2012),
  421--442, \href{http://arxiv.org/abs/1203.3034}{arXiv:1203.3034}.

\bibitem{pilcapaper}
Pilca M., K\"ahlerian twistor spinors, \href{http://dx.doi.org/10.1007/s00209-010-0668-7}{\textit{Math.~Z.}} \textbf{268} (2011),
  223--255, \href{http://arxiv.org/abs/0812.3315}{arXiv:0812.3315}.

\bibitem{SW3}
Seiberg N., Witten E., Monopoles, duality and chiral symmetry breaking in
  {$N=2$} supersymmetric {QCD}, \href{http://dx.doi.org/10.1016/0550-3213(94)90214-3}{\textit{Nuclear Phys.~B}} \textbf{431} (1994),
  484--550, \href{http://arxiv.org/abs/hep-th/9408099}{hep-th/9408099}.

\bibitem{wang}
Wang M.Y., Parallel spinors and parallel forms, \href{http://dx.doi.org/10.1007/BF00137402}{\textit{Ann. Global Anal.
  Geom.}} \textbf{7} (1989), 59--68.

\bibitem{SW2}
Witten E., Monopoles and four-manifolds, \href{http://dx.doi.org/10.4310/MRL.1994.v1.n6.a13}{\textit{Math. Res. Lett.}} \textbf{1}
  (1994), 769--796, \href{http://arxiv.org/abs/hep-th/9411102}{hep-th/9411102}.

\end{thebibliography}
\end{document}